\documentclass{article}

\RequirePackage{amsthm,amsmath,amsfonts,amssymb}

\RequirePackage[numbers]{natbib}

\RequirePackage{graphicx}

\DeclareMathOperator*{\argmin}{argmin}

\begin{document}

\title{\bf Closed-form Tight Bounds and Approximations for the Median of a Gamma Distribution}
\author{Richard F. Lyon \\
  dicklyon@acm.org}
\maketitle

\begin{abstract}
We show how to find upper and lower bounds to the median of a gamma distribution, over the entire range of shape parameter $k > 0$, that are the tightest possible bounds of the form $2^{-1/k} (A + Bk)$, with closed-form parameters $A$ and $B$.  The lower bound of this form that is best at high $k$ stays between 48 and 50 percentile, while the uniquely best upper bound stays between 50 and 55 percentile.  We show how to form even tighter bounds by interpolating between these bounds, yielding closed-form expressions that more tightly bound the median.  Good closed-form approximations between the bounds are also found, including one that is exact at $k = 1$ and stays between 49.97 and 50.03 percentile.
\end{abstract}




\section{Introduction}
\label{sec:intro}

The known Laurent series for the median of a gamma distribution produces some upper and lower bounds, and very accurate approximations at high enough values of the shape parameter $k$, but poor results near $k=1$ and below.  The known approximate and bounding results for $k$ close to zero do not extend to high $k$.  We find new bounds that correspond to percentiles close to 50 across the entire range of shape parameter, and apply them to making even tighter bounds as well as better approximation formulae.  An asymptotic approximation and bound at low $k$, as opposed to the usual focus on high $k$, provides the key functional form needed for this approach.

\section{Problem formulation}
\label{sec:formulation}

The gamma distribution PDF is $\frac{1}{\Gamma(k) \theta^k} x^{k - 1} e^{-\frac{x}{\theta}}$, but we'll use $\theta = 1$ because both the mean and median simply scale with this parameter.  Thus we use this PDF with just the shape parameter $k$, with $k > 0$ and $x \ge 0$:

$$p(x) = \frac{1}{\Gamma(k)} x^{k - 1} e^{-x}$$ 

The mean of this distribution, $\mu$, is well known to be $ \mu = k $.  The median $\nu$ is the value of $x$ at which the CDF equals one-half:

$$\frac{1}{2} = \int_0^\nu p(x) dx = \int_0^\nu \frac{x^{k - 1} }{\Gamma(k)} e^{-x} dx$$

This equation has no easy solution, but the median is well known to be a bit below the mean, bounded by \citep{chen1986bounds}

$$ k - \frac{1}{3} < \nu < k \quad 
\mathrm{and} \quad  0 < \nu $$

Bounds that are tighter in some part of the shape parameter range can be obtained from the known Laurent series partial sums, or from the low-$k$ asymptote and bounds of Berg and Pedersen \cite{berg2006chen}.

We seek upper and lower bounds that are tighter, especially in the middle part of the $k$ range, than are previously known. Further, we seek simple approximation formulae for the median, leveraging these bounds.

\section{Prior work}
\label{sec:prior}

A Laurent series for $\nu$ with rational coefficients has been discovered, with deep connections to some math by Ramanujan.  Choi \cite{choi1994medians} applied Ramanujan's work to this particular question, providing 4 coefficients (through the $k^{-3}$ term). Berg and Pedersen \cite{berg2006chen}, based on work by Marsaglia \cite{marsaglia1986c249}, extended this to 10 coefficients.  Neither commented on the radius of convergence, which appears to be in the neighborhood of $k = 1$.  So for large enough $k$ and $N$ the series yields excellent approximations, but for $k < 1$ it is useless. 
$$ \nu \approx k + \sum_{j = 0}^N a_j k^{-j}$$
with $ a_j = \{ 
\frac{-1}{3}, 
\frac{2^3}{3^4 \cdot 5}, 
\frac{2^3 \cdot 23}{3^6 \cdot 5 \cdot 7},
\frac{2^3 \cdot 281}{3^9 \cdot 5^2 \cdot 7},  
\frac{-2^3 \cdot 17 \cdot 139753}{3^{13} \cdot 5^3 \cdot 7 \cdot 11}, 
\frac{-2^3 \cdot 708494947}{3^{15} \cdot 5^3 \cdot 7^2 \cdot 11 \cdot 13}
\dots \} $.  Thus, where Choi \cite{choi1994medians} had 144 instead of the correct $2^3 \cdot 23 = 184$:
$$\nu = k - \frac{1}{3} + \frac{8}{405 k} + \frac{184}{25515 k^2} + O\left(\frac{1}{k^3}\right)  $$

Partial sums of the Laurent series are not generally bounds, but the first two ($k$ and $k - 1/3$) are upper and lower bounds, respectively, and the sums ending with $-3$ and $-5$ powers of $k$ are also upper and lower bounds, respectively.  

Berg and Pedersen \cite{berg2006chen} also derived an asymptote for small $k$, which we call $\nu_0$:
$$\nu \approx \nu_0 = e^{-\gamma}\,2^{-1/k}$$
where $\gamma \approx 0.577216$ is the Euler--Mascheroni constant.  This asymptote is a lower bound, as we will show.  Their factor $2^{-1/k}$ is the key to good approximations and bounds, and we divide by it to reduce the dynamic range in some of our later plots.  For the range $0.01 < k < 100$, this reduces the dynamic range we need to work with by nearly 30 orders of magnitude---or 300 orders of magnitude for $k$ down to 0.001---but we still need log--log plots to show the bounds and approximations, and their errors, across these ranges.

Berg and Pedersen \cite{berg2006chen} also provide an upper bound $\nu < e^{-1/3k}k$ that is just above $k - \frac{1}{3}$, and a lower bound $\nu > 2^{-1/k}k$.

\begin{figure}
  \centering\includegraphics[width=3.85in]{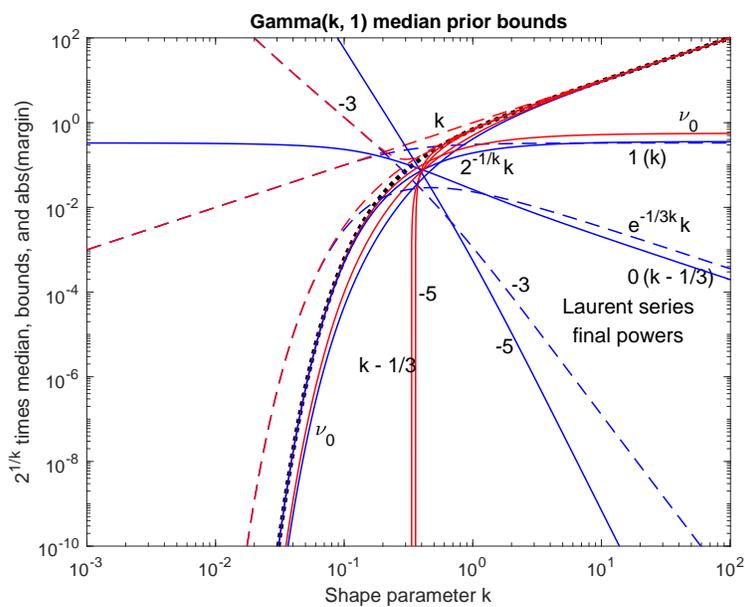}
  \caption{Previously published bounds (red; lower bounds solid, upper bounds dashed) for the median of a gamma distribution (black dotted), are good at high $k$ or low $k$, but not both.  Their margins (errors) are shown in blue (corresponding line styles).  At the left, at $k = 0.01$, the median is near $10^{-30}$, and at $k = 0.001$ near $10^{-300}$.}
  \label{fig:priorart}
\end{figure}

These previous known bounds are illustrated in Figure \ref{fig:priorart}, where upper and lower bounds and their errors (or ``margins'') are distinguished by different line styles.

Others have shown good bounds and approximations where $k$ is an integer, that is, for the Erlang distribution \citep{adell2008ramanujan, you2017approximation, chen2017median}.  These results do not extend to the  low-$k$ region.

The approach of approximating functions by interpolating between upper and lower bounds has been discussed by Barry \cite{barry2000approximation}, who used minimax optimization to find numeric parameters in an interpolation-between-bounds approximation to the exponential integral.  We are not aware of an interpolation approach being used to find improved closed-form bounds.

\section{Tight upper and lower bounds}
\label{sec:approach}

The Berg lower bounds $2^{-1/k}e^{-\gamma}$ and $2^{-1/k}k$ are tight within their families $2^{-1/k}A$ and $2^{-1/k}Bk$, but not tight within the wider two-parameter family $2^{-1/k}(A + Bk)$.  We show that the sum of these two is the uniquely tight upper bound in that family, and that there is a range of tight lower bounds in the family.

To improve on the lower bounds, and to motivate the family that we consider further, first solve for $\nu$ in a simple approximation to the distribution's integral, using $e^{-x} < 1$ for $x > 0$:
$$\frac{1}{2} < \int_0^\nu \frac{x^{k - 1} }{\Gamma(k)} dx$$
$$ \nu > 2^{-1/k}\ \Gamma(k+1)^{1/k} $$
which is a tight lower bound, and is a good approximation for $k < 0.1$, but not so great at high $k$---and not what we consider a closed form, due to the gamma function.  This expression resembles the ``quantile mechanics'' boundary condition for the gamma distribution from Steinbrecher and Shaw \cite{steinbrecher2008quantile}, and converges with the Berg asymptote at low $k$, but we can improve Berg's result by utilizing another term.  A symbolic calculus system finds for us the next Taylor series terms about $k = 0$ for the power of the gamma function:
$$ \Gamma(k+1)^{1/k} \approx  e^{-\gamma} + \frac{e^{-\gamma}\pi^2}{12} k - 0.035 k^2 $$

\begin{figure}
  \centering\includegraphics[width=3.85in]{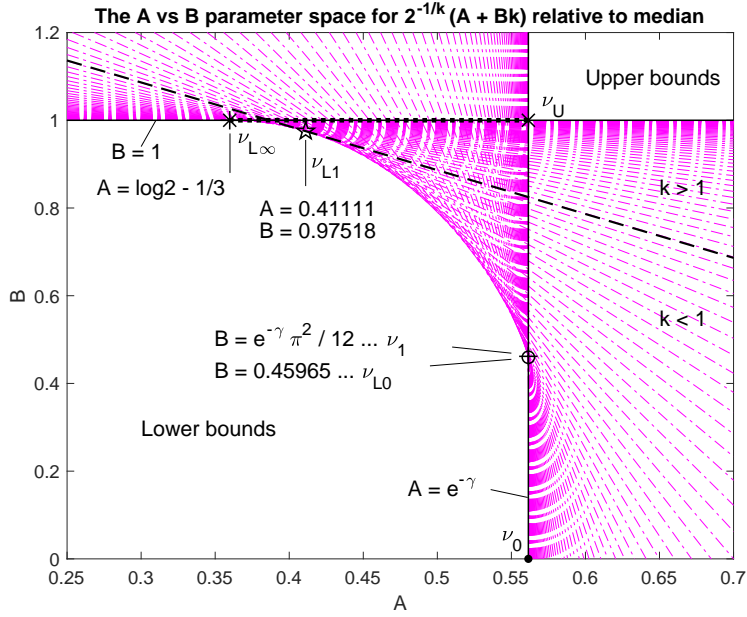}
  \caption{The $A$--$B$ parameter space is shaded with dash-dot lines where $\nu(k) 2^{1/k} = A + Bk$, for a set of very small to very large $k$ values in geometric progression (using numerically computed $\nu(k)$ values). Key values of $A$ and $B$ are indicated. Points outside (or on the edge) of the shaded region represent bounds, while points inside the shaded region represent functions that cross the median function.  There is an obvious uniquely tight upper bound $\nu_U$, and a curved locus of tight lower bounds from $\nu_{L0}$, which is tightest near $k=0$, to $\nu_{L\infty}$, which is tightest for high $k$.  One point (pentagram) on the curved locus represents a lower bound $\nu_{L1}$ that is tight at $k=1$, for which $A + B = 2 \log 2$ (which is the equation of the dashed line).  The point $\nu_1$ represents a good asymptotic approximation close to $\nu_{L0}$, but not a bound; see the next figure.  The dotted line from $\nu_U$ to $\nu_{L\infty}$ at $B = 1$ intersects the lines for all $k$ in monotonic order, with $A$ decreasing while $k$ increases.}
  \label{fig:space}
\end{figure}

So we have this improved approximation, which has much less relative error than Berg's at low $k$, and has a high-$k$ behavior nearly proportional to $k$ (but is no longer a lower bound because we made it larger by ignoring a next negative term):
$$ \nu_1 = 2^{-1/k} \left(  e^{-\gamma} + \frac{e^{-\gamma}\pi^2}{12} k \right) $$

Inspired by this asymptotic approximation, we consider members of this family of functions,  with coefficients $A$ and $B$, and analyze which ones are bounds:
$$\tilde{\nu} = 2^{-1/k}\,(A + B\,k)$$

A graphical characterization of this family is most informative.  Given some values of $k$ and corresponding numerical $\nu(k)$, we can find the lines $A + Bk = \nu(k) 2^{1/k}$ in $A$--$B$ space, and plot them---see Figure \ref{fig:space}.  Regions full of lines are not bounds, and regions without lines are where bounds are found (including some of Berg's bounds); we're interested in the boundaries between these regions, where tight bounds are to be found.

\begin{figure}
  \centering\includegraphics[width=3.85in]{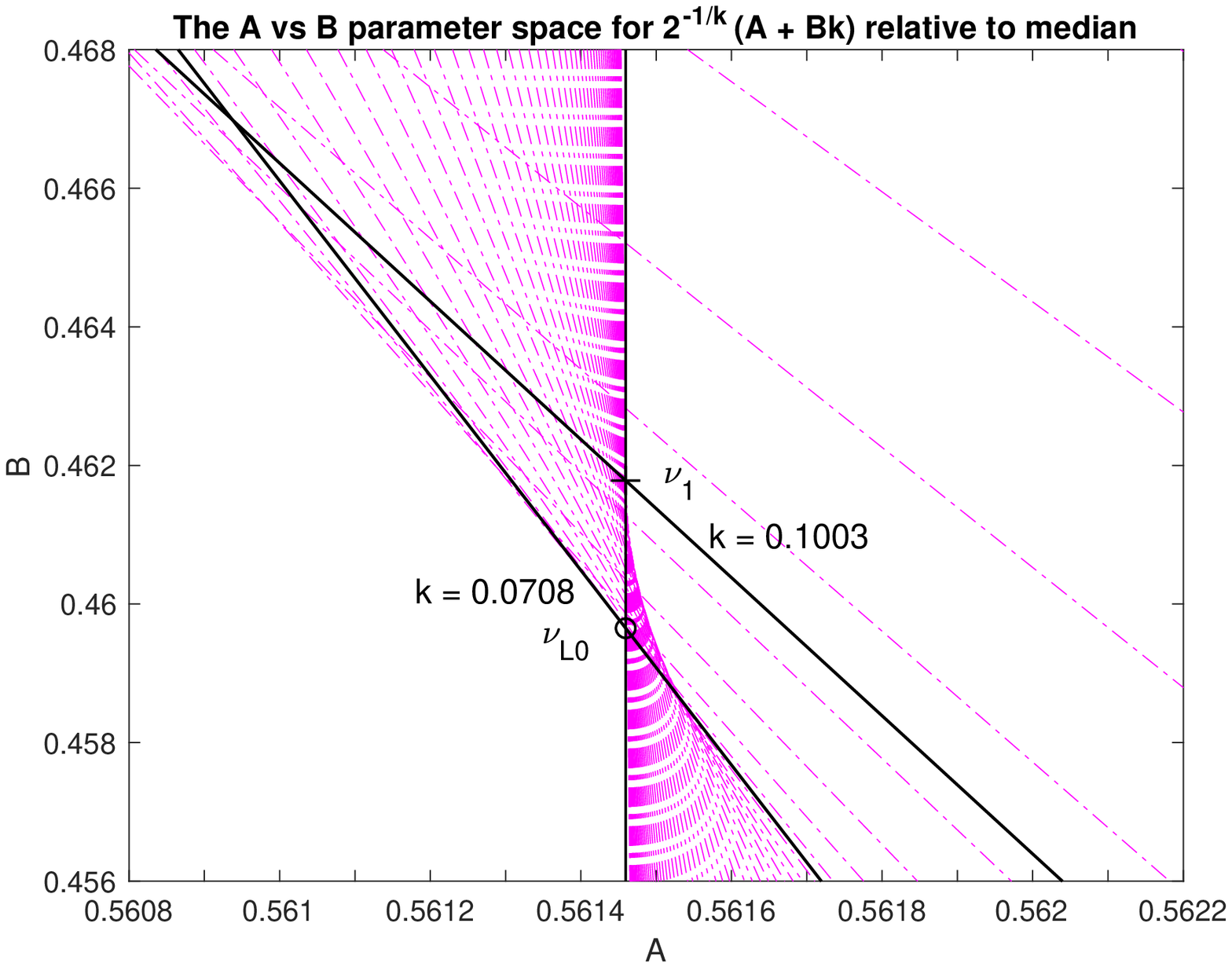}
  \caption{Zooming in to $\nu_1$ and $\nu_{L0}$, note that the point  with $B = e^{-\gamma}\pi^2/12$, which we got from the Taylor series of the power of the gamma function, is actually inside the shaded area, so does not represent a bound; but a point at slightly lower $B = 0.45965$ is on the edge, so represents a lower bound.  These points give zero error at approximately $k = 0.1003$ and $k = 0.0708$, respectively (see the next figure).  We do not have analytic formulations for these numeric and graphical observations.}
  \label{fig:zoom}
\end{figure}

\begin{figure}
  \centering\includegraphics[width=3.85in]{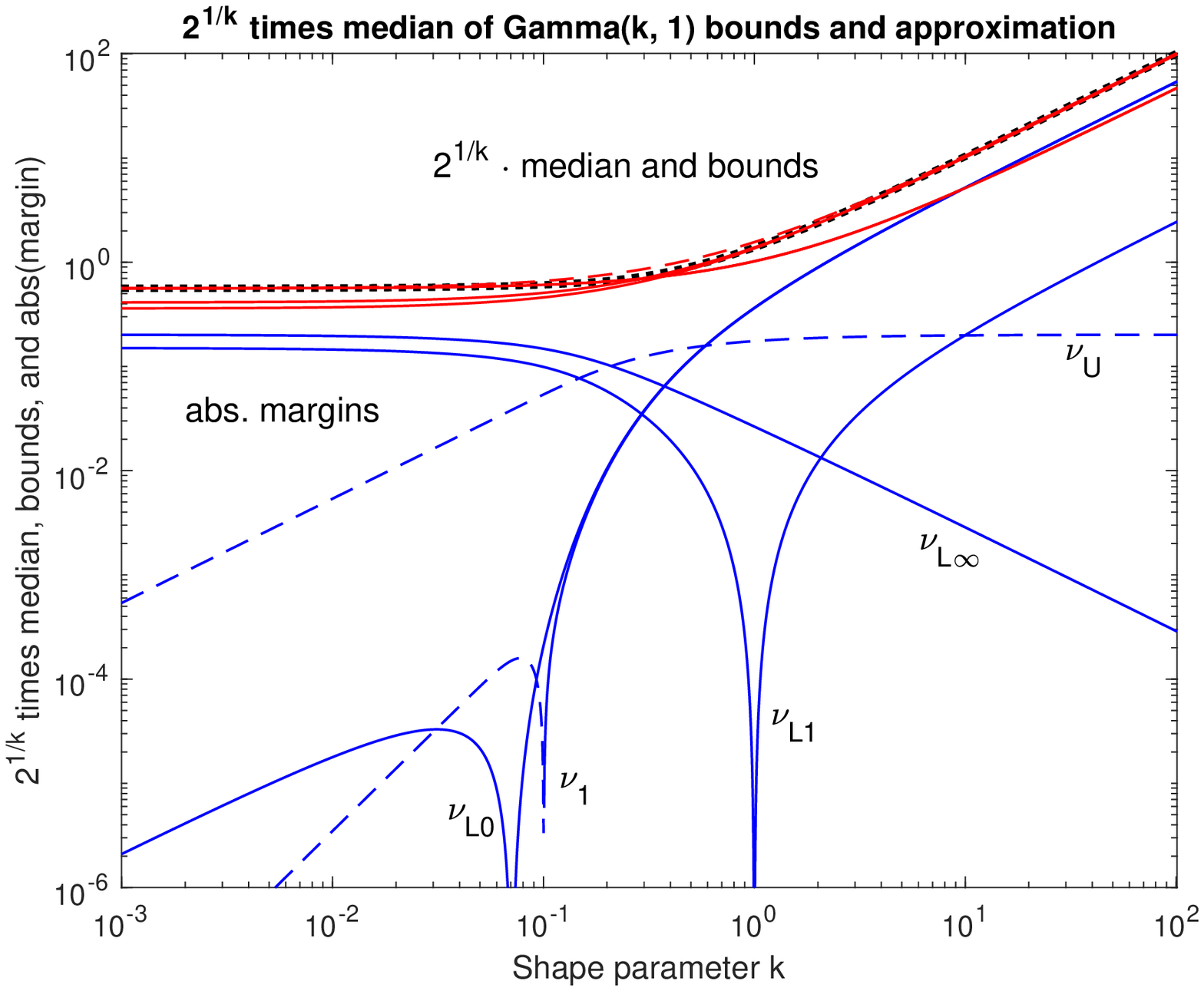}
  \caption{The lower bounds $\nu_{L0}$, $\nu_{L1}$, and $\nu_{L\infty}$ (solid), and upper bound $\nu_U$ (dashed) are shown in red over the ideal median (black heavy dots), with their absolute errors in blue, all premultiplied by $2^{1/k}$ to reduce the required plot range.  The approximation $\nu_1$, which is not a bound, is also shown; note that its error curve changes from solid to dashed at the cusp, while the errors for $\nu_{L0}$ and $\nu_{L1}$ have $\log(0)$ cusps where the error grazes zero but does not change sign. The $k$ parameters at these cusps correspond to the sloped lines indicated in the previous figures.}
  \label{fig:components}
\end{figure}

The improved asymptote $\nu_1$ is great at low $k$, but is neither an upper nor a lower bound, as shown in \ref{fig:zoom}.  We can modify it to approach $k - \frac{1}{3}$ at high $k$ by a few adjustments, via this asymptotic approximation that we get from a symbolic calculus system:
$$k\,2^{-1/k} = k - \log 2 + O(k^{-1}) $$

Thus we find this approximation for high $k$, which is a lower bound as illustrated in Figure \ref{fig:space}:
$$\nu_{L\infty} = 2^{-1/k} \,\left(\log 2 - \frac{1}{3} + k \right)$$

A compromise approximation for low $k$ mixes these two, differing from the high-$k$ approximation in only the $A$ coefficient, leaving a result consistent to the same order as Berg's at low $k$, and forming an upper bound as illustrated in \ref{fig:space}:
$$\nu_U = 2^{-1/k} \left(  e^{-\gamma} + k \right)$$

This mixed approximation has absolute and relative errors approaching zero at low $k$, and relative error approaching zero at high $k$; but the absolute error remains high, near $\log 2 - \frac{1}{3} - e^\gamma \approx 0.20$, at high $k$. These approximations and their errors are illustrated in Figure \ref{fig:components}.

To support the graphical/numerical observation that $\nu_{L\infty}$ and $\nu_U$ are lower and upper bounds, respectively, of the true median ($ \nu_{L\infty} < \nu < \nu_U $), we examine their asymptotic behaviors in more detail. At low $k$, it is easy to see, using $\log 2 - \frac{1}{3} \approx 0.359814 < e^{-\gamma} \approx 0.561459$, and $\frac{e^{-\gamma}\pi^2}{12} \approx 0.461781 < 1$, that these differences are positive, for $k \rightarrow 0$:

$$2^{1/k} \left(\nu - \nu_{L\infty} \right) = e^{-\gamma} - \left(\log 2 - \frac{1}{3}\right) + O(k)> 0$$
$$2^{1/k} \left( \nu_U - \nu \right) =  k \left( 1 -  \frac{e^{-\gamma}\pi^2}{12} \right)   + O(k^2)> 0$$ 

At high $k$, a symbolic calculus system gives us for $\nu_{L\infty}$:
$$2^{-1/k}k = k - \log 2 + \frac{\log^2 2}{2}\,k^{-1} + O(k^{-2})$$
which we can use to construct comparisons to the Laurent series terms for $\nu$ \citep{choi1994medians,berg2006chen}.  Again we find positive differences, with $\frac{\log 2}{3} - \frac{\log^2 2}{2} \approx -0.009177 < \frac{8}{405}$, and $\log 2 - e^{-\gamma} \approx 0.131688 < \frac{1}{3}$, for $k \rightarrow +\infty$:

$$ \nu - \nu_{L\infty} = \left(\frac{8}{405}  -  \left( \frac{\log 2}{3} - \frac{\log^2 2}{2} \right)  \right) \,k^{-1} + O(k^{-2}) > 0 $$
$$ \nu_U -\nu =   e^{-\gamma} - \log 2  + \frac{1}{3} + O(k^{-1}) > 0$$

In addition to these high-$k$ and low-$k$ asymptotic results, we can show the inequalities also hold at $k = 1$ where $\frac{1}{2}\left(\log 2 - \frac{1}{3} + 1 \right) < \log 2 < \frac{1}{2} (e^{-\gamma} + 1)$, but otherwise we're relying on the graphical and numerical results. Since there is considerable margin in the asymptotes, and the median is well behaved (unique, monotonic, and smooth, with positive second derivative for all $k > 0$ \citep{berg2006chen,berg2008convexity}), this seems reliable enough in concluding that these are bounds (but mathematicians are invited to interpret these as conjectures to be proved).
The positivity constraints would not all hold, since higher-order terms would not cancel, if $A_{L\infty}$ or $B_{L\infty}$ were any higher, or if $A_U$ or $B_U$ were any lower.  In that sense, these upper and lower bounds are proved tight.

For the lower bound $\nu_{L1}$ that is tight at $k = 1$, both the value and the slope need to match the true median.  The value is the median of the exponential distribution, $\nu(1) = \log 2$.  The slope $v^{\prime }(k)$, which is somewhat more troublesome to work out, but is tractable at the special point $k = 1$, is:
$$v^{\prime }(k)|_{k=1} = \gamma - 2 \textrm{Ei}(-\log 2) - \log \log 2 \approx 0.9680448$$
which is a mathematical expression with a definite value, but is not a closed form due to the exponential integral, so still requires a numerical approach to evaluate it.  Therefore, we describe this bound with approximate numerical parameters instead of closed-form analytic expressions.

$$ B = 2 \left( v^{\prime }(k)|_{k=1}  - \log^2 2 \right) \approx 0.9751836$$
$$ A = 2 \log 2 - B \approx 0.4111107$$

The lower bound $\nu_{L0}$ is in worse shape, as we have to search for the $k$ value that gives the lowest $B$ value with $A = e^{-\gamma}$.  So its $B$ parameter has no concise mathematical expression, but can be computed to high precision: $B \approx 0.4596507$.

\begin{table}[t]
 \caption{Comparison of several median bounds and asymptotes.}
  \centering
  \begin{tabular}{lccll}
    Version & \multicolumn{2}{c}{ } & Description  \\ [3 pt]
    \hline
    $\nu$ & & & True median of gamma distribution \\ [3 pt]
    \multicolumn{2}{l}{$e^{-1/3k}k$} & & Berg's upper bound, high-$k$ asymptote \\ [3 pt]
    \multicolumn{2}{l}{$2^{-1/k}\Gamma(k+1)^{1/k}$} & & New lower bound, low-$k$ asymptote \\ [3 pt]
    \hline
    & $A$ & $B$ & parameters for the form $2^{-1/k}(A + Bk)$ \\ [3 pt]
    \hline
    $2^{-1/k}k$ & $0$ & 1 & Berg's lower bound \\ [3 pt]
    $\nu_0$ & $e^{-\gamma}$ & 0 & Berg's asymptote, a lower bound \\ [3 pt]
    $\nu_1$ & $e^{-\gamma}$ & $\frac{e^{-\gamma}\pi^2}{12}$ & Improved low-$k$ asymptote; not a bound \\ [3 pt]
    $\nu_{L0}$ & $e^{-\gamma}$ & $0.4596507$ & New tight lower bound, best at low $k$  \\ [3 pt]
    $\nu_{L1}$ & $0.4111107$ & $0.9751836$ & New tight lower bound, tangent at $k = 1$  \\ [3 pt]
    $\nu_{L\infty}$ & $\log 2 - \frac{1}{3}$ & $1$ & New tight lower bound, best at high $k$ \\ [3 pt]
    $\nu_U$ & $e^{-\gamma}$ & $1$ & New uniquely tight upper bound \\ [3 pt]
  \hline

  \end{tabular}
  \label{tab:compareAB}
\end{table}

The bounds and asymptotic approximations discussed here are summarized in Table \ref{tab:compareAB}.  In subsequent sections we focus primarily on the new upper and lower bounds with closed-form coefficients, $\nu_U$ and $\nu_{L\infty}$, as a basis for even tighter closed-form bounds.

\begin{figure}
  \centering\includegraphics[width=3.85in]{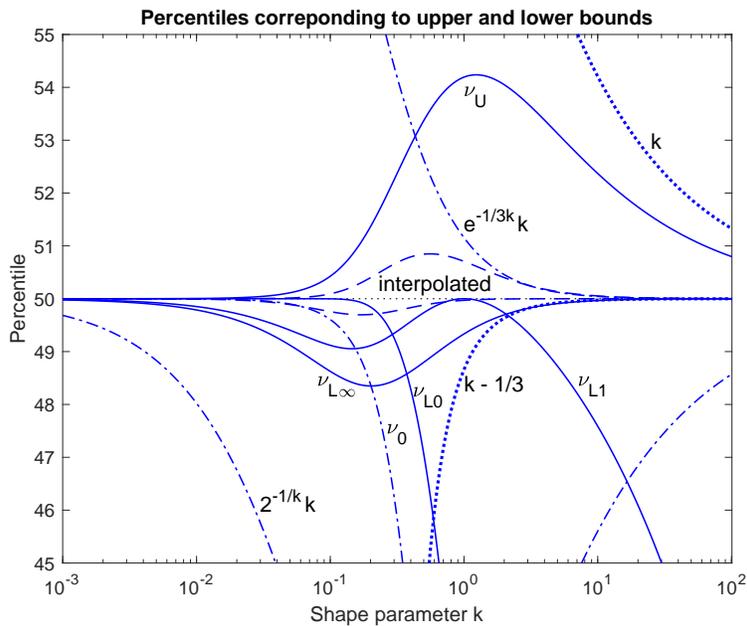}
  \caption{The percentiles achieved by four new median bounds of the form $2^{-1/k}(A + Bk)$ (solid curves) are plotted, along with the linear bounds $k$ and $k - \frac{1}{3}$ (dotted), upper and lower bounds from Berg and Pedersen \cite{berg2006chen} (dash-dot), and a pair of closer bounds formed by interpolation between $\nu_U$ and $\nu_{L\infty}$ using a one-parameter rational function (dashed).  The bounds $2^{1/k}k$, $\nu_U$, $\nu_{L\infty}$, and the interpolated bounds converge on 50th percentile at both low and high $k$, while the other six do not.  Both the upper and lower interpolated bounds are close to $\nu_U$ at low $k$ and close to $\nu_{L\infty}$ at high $k$; tighter such interpolated bounds, developed in a later section, would crowd the center of the graph.}
  \label{fig:percentiles}
\end{figure}

Figure \ref{fig:percentiles} shows the percentile values achieved by these bounds and approximations, compared to the ideal 50\% that defines the median---$\nu_{L\infty}$ always comes in between 48\% and 50\% and $\nu_U$ always between 50\% and 55\% (percentiles are calculated in Matlab as \mbox{100*gammainc(x, k)}, using the normalized lower incomplete gamma function that is the CDF for our PDF).

The coefficient of $k^{-1}$ for $\nu_{L\infty}$ is negative, so the lower bound $\nu_{L\infty}$ slightly violates the $k - \frac{1}{3}$ lower bound in spite of having asymptotically zero absolute error.  That is, for $k > 3.021$, it's a looser lower bound and a worse approximation than $k - 1/3$, even though it is the tightest lower bound of the form we're considering.  On the other hand, $\nu_U$ is a much tighter upper bound than $k$ is.

\section{Formulae for tighter bounds}
\label{sec:tighter}

\begin{figure}
  \centering\includegraphics[width=3.85in]{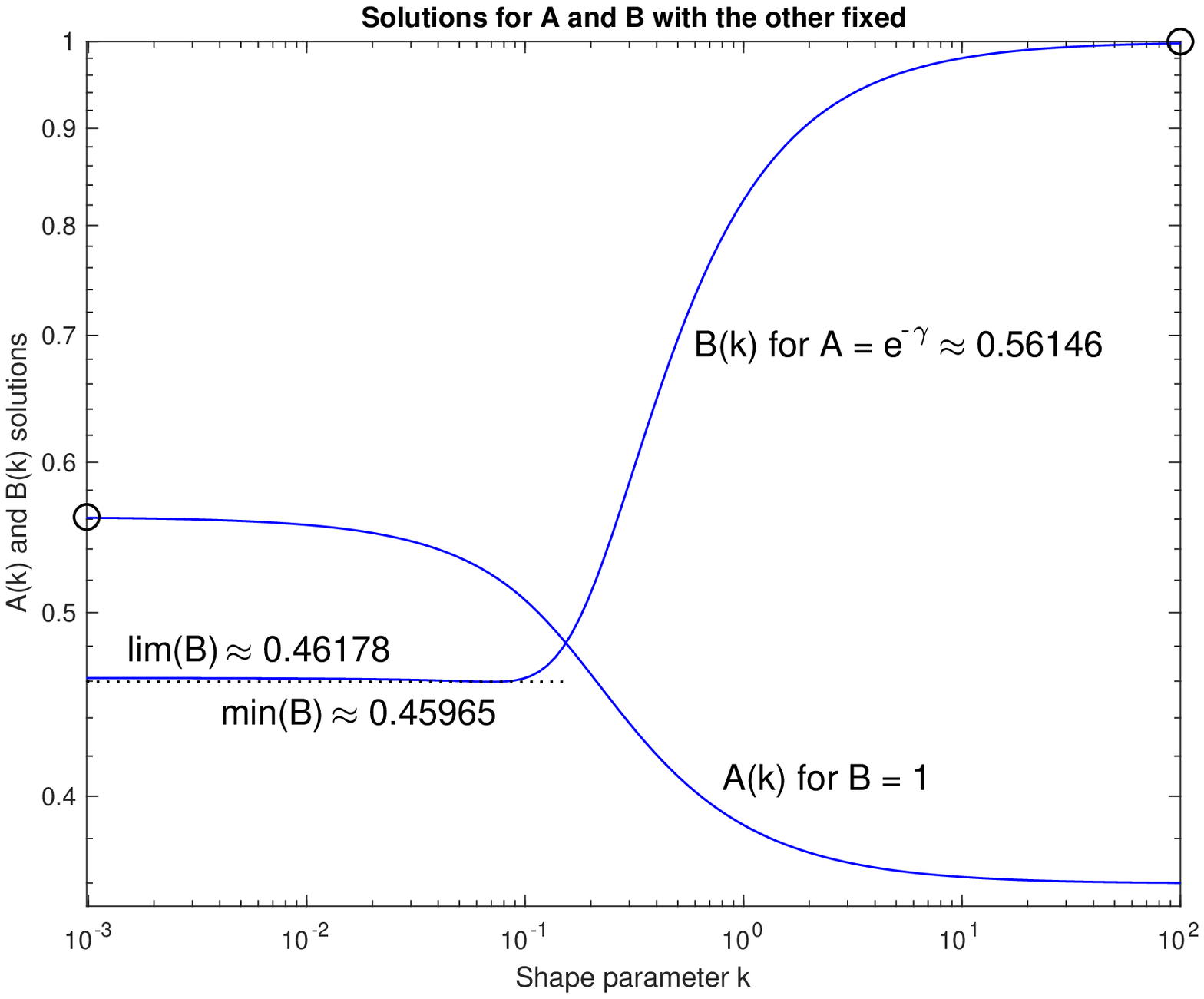}
  \caption{Functions $A(k)$ and $B(k)$, either of which can solve $\nu = 2^{-1/k}(A + B k)$, with the other constant at the limiting values indicated by the circles, the parameters of $\nu_U$; modeling these curves can lead to better bounds or approximations.}
  \label{fig:A_B_functions}
\end{figure}

Letting $A$ and $B$ be functions of $k$, rather than constants, allows tighter bounding expressions (and potentially exact expressions) for $\nu(k)$, but not enough structure.  Allowing only one of them to vary, and tying the other to values used in the tight bounds above, allows a more constrained space of bounds.

With $B_{L\infty} = B_U = 1$, we can express the median exactly as $\nu = 2^{-1/k}(A(k) + k)$, for some smooth positive real function $A(k)$ that runs from a limit of $A_U = e^{-\gamma}$ as $k \rightarrow  0$ to $A_{L\infty} = \log 2 - \frac{1}{3}$ as $k \rightarrow +\infty$; it is apparently monotonic.  

Alternatively, using $A_{L0} = A_U = e^{-\gamma}$, the formula $\nu = 2^{-1/k}(e^{-\gamma} + B(k)k)$ has a smooth positive but non-monotonic function $B(k)$ that runs between limits $\frac{e^{-\gamma} \pi^2}{12}$ and $1$, but drops a little below its low-$k$ limit before increasing.  

This approach converts the problem of finding tighter bounds to the median to the problem of finding closed-form expressions to bound these well-behaved functions.  Calculating $A(k)$ and $B(k)$ numerically to high precision is easy when the median is calculated; see Figure \ref{fig:A_B_functions}.  For the rest of this paper, we focus on $A(k)$, since it is monotonic and more nearly symmetric on a $\log k$ axis, and because it corresponds to interpolation between closed-form bounds.

\section{Interpolators}
\label{sec:interps}

The function $A(k)$ introduced above can be represented in terms of an interpolation function $g(k)$ that runs monotonically from 0 at low $k$ to 1 at high $k$:
$$A(k) = g(k)A_{L\infty} + (1 - g(k))A_U = A_U - g(k)(A_U - A_{L\infty}) $$
$$ = e^{-\gamma} - g(k) \left(e^{-\gamma} - \log 2 + \frac{1}{3} \right)  $$
And $g(k)$ is therefore also the function that interpolates between the bounds, allowing us to write the median in these convenient ways:
$$\nu = g(k)\nu_{L\infty} + (1 - g(k))\nu_U$$
$$\nu = 
2^{-1/k} \left( e^{-\gamma} - g(k)  \left(e^{-\gamma} - \log 2 + \frac{1}{3} \right)  + k \right)$$

The ideal interpolator can be found numerically:
$$ g(k) = \frac{A_U - A(k)}{A_U - A_{L\infty}} = \frac{\nu_U - \nu}{\nu_U - \nu_{L\infty}} $$

It can be interesting to bound or otherwise approximate $g(k)$ with various $\tilde{g}(k)$.
In approximating the ideal with an interpolator $\tilde{g}(k)$, we achieve absolute and relative error of the median estimate approaching zero at low $k$ if $\tilde{g}(k) = 0 + O(k)$, and at high $k$ if $\tilde{g}(k) = 1 - O(k^{-1})$.  But we might want to do better, matching the asymptotic slopes of the ideal interpolator to match the median to a higher order; or we might want to match the exact known value $\nu = \log 2$ at $k = 1$.  So we analyze these properties of the ideal interpolator, and give them names.  At low $k$:
$$ P_0 = \frac{dg}{dk}  = \frac{1 - \frac{e^{-\gamma} \pi^2}{12}}{e^{-\gamma} - \log 2 + \frac{1}{3}} \approx 2.66913 $$
At high $k$:
$$ P_\infty = -\frac{dg}{d\frac{1}{k}} = \frac{\frac{8}{405} + e^{-\gamma} \log 2 - \frac{\log^2 2}{2}}{e^{-\gamma} - \log 2 + \frac{1}{3}} - \log 2 \approx 0.143472$$
And at $k = 1$:
$$ P_1 = g(1) = \frac{1 + e^{-\gamma} - 2\log 2}{e^{-\gamma} - \log 2 + \frac{1}{3}} \approx 0.868678 $$

\begin{figure}
  \centering\includegraphics[width=3.85in]{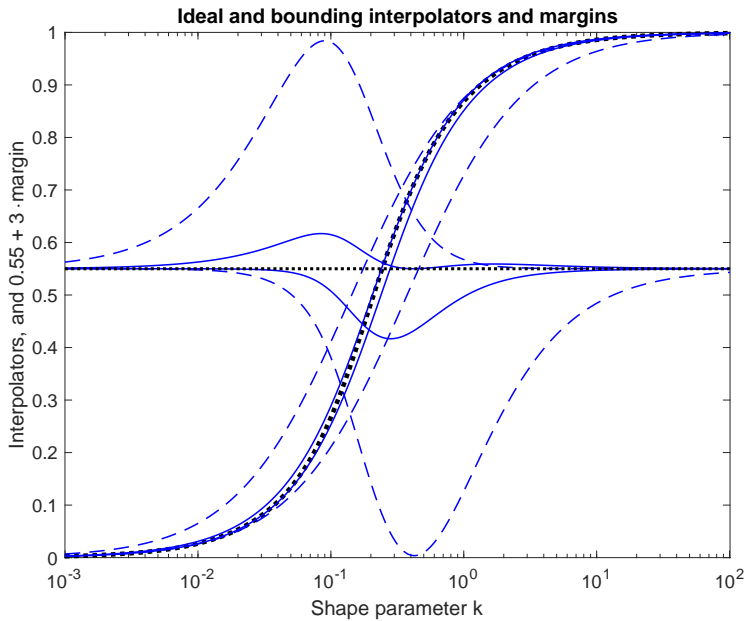}
  \caption{The ideal interpolator $g(k)$ (heavy dotted sigmoid) is compared with upper and lower bounds $\tilde{g}(k)$; their margins $\tilde{g}(k) - g(k)$ are also plotted, magnified and displaced, with the same curve styles.  The curves with largest absolute margins (dashed), which correspond to the interpolated bounds shown in Figure \ref{fig:percentiles}, are for the first-order rational-function interpolator $\frac{k}{k + b_0}$, while the curves with smaller margins (solid) are for arctan interpolators $\frac{2}{\pi}\tan^{-1} \frac{k}{b}$.  In each case, the one parameter ($b_0$ or $b$) is chosen to give a tight bound (analytically in closed form in three of the four cases). 
  Lower bounds of $g(k)$ make upper bounds of $\nu(k)$, and vice versa.}
  \label{fig:interpbounds}
\end{figure}

\begin{figure}
  \centering\includegraphics[width=3.85in]{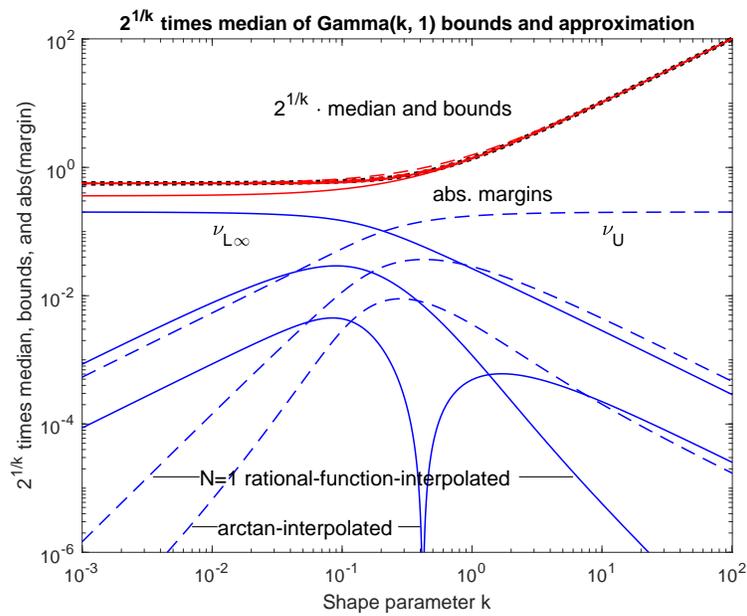}
  \caption{The absolute margins of the interpolated bounds are smaller than those of the bounds they started from.  Compare Figure \ref{fig:components}. The generally smallest margins are for the arctan interpolators, and the intermediate for the $N=1$ rational-function interpolators.}
  \label{fig:interp_margins}
\end{figure}

How such approximation goals relate to bounds is not immediately clear.  We construct some examples.
Figure \ref{fig:interpbounds} shows the ideal interpolator, computed numerically, and compares it to bounding interpolators of the forms $\tilde{g}(k) = \frac{k}{k + b_0}$ and $\tilde{g}(k) = \frac{2}{\pi}\tan^{-1} \frac{k}{b}$, with parameters $b_0$ and $b$ chosen to yield tight upper and lower bounds, as discussed in the next two sections.  The effects of these interpolator bounds on the median bounds is shown in Figure \ref{fig:interp_margins}.

\section{Rational-function interpolators}

Consider rational functions as interpolators, of the form
$$\tilde{g}_N(k) = \frac{\sum_{n = 1}^{N-1} a_i k^i + k^N}{\sum_{n = 0}^{N-1} b_i k^i + k^N}$$

For $N=1$, the only parameter is $b_0$, so we have a one-parameter family:
$$\tilde{g}_1(k) = \frac{k}{b_0 + k}$$

For $N=2$ we have three parameters:
$$g_2(k) = \frac{a_1k + k^2}{b_0 + b_1k + k^2}$$
and so forth.

We can easily constrain the coefficients to match the properties of the ideal interpolator.  
At low $k$:
$$ \frac{a_1}{b_0}  = P_0$$
At high $k$
$$b_{N-1} - a_{N-1} = P_\infty$$
And at $k = 1$:
$$ \frac{\sum_{n = 1}^{N-1} a_i + 1}{\sum_{n = 0}^{N-1} b_i + 1} = P_1$$

For $N = 1$, the low asymptote is tightly approached with $b_0 = \frac{1}{P_0} $, yielding an upper bound for the median (lower bound for the interpolator $g_{\mathrm{ideal}}$).  Or the high asymptote is tightly approached with $b_0 = P_\infty $, yielding a lower bound for the median (upper bound for $g(k)$).  These bounds are illustrated in Figure \ref{fig:interpbounds}.
For $b_0$ between these values, the resulting interpolated function is not a bound, but is exact at one value of $k$, for example at $k = 1$ with $b_0 = \frac{1}{P_1} - 1 $.

For $N = 2$, there are enough parameters to use any or all of the three constraints, but it's not immediately clear which sets of constraints can lead to bounds.  Certainly using all three constraints does not lead to a bound, but to an interesting approximation.  As we did with the $A$ versus $B$ space, we can investigate the locus of parameter solutions for each $k$, and examine the edges of these locus-filled areas for tight bounds.  Reducing the space to 2D by constraining one asymptote or the other allows a graphical approach, but the results are not better than a one-parameter arctan interpolator.

For $N \geq 3$, excellent approximations and bounds are possible, but with so many parameters are not very interesting.

\section{Arctan interpolators}
\label{sec:arctan}

Getting a tighter bound or better approximation from the rational-function family requires at least a handful of parameters.  An alternative approach is to find a one-parameter shape that fits better.  We have found that the arctan shape with parameter $b$ (like $b_0$ in the one-parameter rational-function interpolator, corresponding to the $k$ value at the midpoint, $\tilde{g}(b) = 0.5$) does a good job:
$$ \tilde{g}_a(k) = \frac{2}{\pi} \tan^{-1} \frac{k}{b} $$

As Figure \ref{fig:interpbounds} shows, the shapes of the arctan interpolators are imperfect, but are much better than the first-order rational function, fitting better in some regions than in others.  Note that both the $N = 1$ rational function and the arctan interpolators are symmetric about their centers (with $\log k$ as the independent variable); they are logistic sigmoid and Gudermannian shapes, respectively.  The ideal that they are to bound or approximate, however, is not quite symmetric in $\log k$.  So a family of not-quite-symmetric interpolators can perhaps do better.

\begin{figure}
  \centering\includegraphics[width=3.85in]{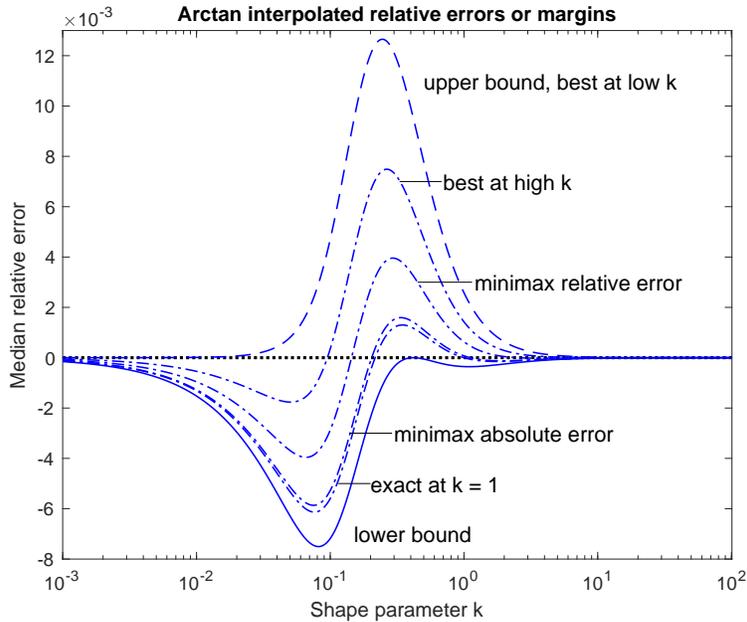}
  \caption{Relative errors of arctan-interpolated approximations (dash-dot curves) between the upper (dashed) and lower (solid) bounds.  These are among the possibly interesting approximations suggested by the tight-bounds approach.  The non-bounding approximations, optimized for different criteria, all have maximum relative errors below 1\%.  }
  \label{fig:interprel}
\end{figure}

Several special $b$ parameter values for the arctan interpolator can be derived to match each of the ideal properties mentioned above.
We can constrain the approximation to pass through the known value $\nu = \log 2$ at $k = 1$, with $b = \cot \left( \frac{\pi}{2} P_1 \right) $, but that does not give a bound.
Or we match the true median at high $k$ to within $O(k^{-2})$ with $b = 
\frac{\pi}{2} P_\infty $.  That also does not give a bound.
Or we can match at low $k$ to within $O(2^{-1/k}\,k^2)$, like $\tilde{\nu}_1$ does, with $b = \frac{2}{\pi} P_0^{-1}$, which yields a lower bound to $g$.

To find an upper bound to $g$ (lower bound to $\nu$), we decrease $b$ until the margin is nonnegative for all $k$, which is at about $b = 0.205282$.  Finding an analytic formulation for that tight bound is a challenge left to others.

Several interpolated bounds and approximations are summarized in Table \ref{tab:compareInterps}.  Most are closed-form analytic expressions; see Figure \ref{fig:interprel} for the relative errors of various arctan interpolator versions.

\begin{table}[t]
 \caption{Comparison of several one-parameter interpolated bounds and approximations $\tilde{g}(k)\nu_{L\infty} + (1 - \tilde{g}(k))\nu_U$, some of which have closed-form parameters.  Bounds are indicated by U or L.}
  \centering
  \begin{tabular}{lccc}
    Version & \multicolumn{3}{c}{Parameter}  \\ [2 pt]
    \hline
    $\tilde{g}_1(k) = \frac{k}{k+b_0}$ & Symbolic $b_0$  &  Numeric $b_0$ \\[5 pt]
    
    best at low $k$ &  $ \frac{e^{-\gamma} - \log 2 + \frac{1}{3}}{1 - \frac{e^{-\gamma} \pi^2}{12}} $  &  0.374654 & U \\[7 pt]

    exact at $k = 1$    
    &  {
    $ \frac{e^{-\gamma} - \log 2 + \frac{1}{3}}{1 + e^{-\gamma} - 2\log 2} - 1 $} 
  & 0.151175 & -- \\[7 pt]

    best at high $k$ &  $ \frac{\frac{8}{405} + e^{-\gamma} \log 2 - \frac{\log^2 2}{2}}{e^{-\gamma} - \log 2 + \frac{1}{3}} - \log 2 $ & 0.143472 & L \\[7 pt]
    
  \hline
    $\tilde{g}_a(k) = \frac{2}{\pi} \tan^{-1} \frac{k}{b} $ & Symbolic $b$  &  Numeric $b$ \\[5 pt]
    
    best at low $k$  
    &  $\frac{24}{\pi} \left( \frac{e^{-\gamma} - \log 2 + \frac{1}{3}}{12 - e^{-\gamma}\pi^2} \right)$  &  
    0.238512 & U \\[7 pt]

    best at high $k$    
    &  {
    $ \frac{\pi}{2} \left( 
  \frac{ \frac{8}{405} + e^{-\gamma} \log 2 - \frac{\log^2 2}{2}} {e^{-\gamma} - \log 2 + \frac{1}{3}} - \log 2 \right) $} 
  & 0.225366 & -- \\[7 pt]

    minimax relative error    
    &  $ \argmin_r \max  \left| \frac{\nu - \tilde{\nu}}{\nu} \right| $
  & 0.21639 & -- \\[7 pt]

    minimax absolute error    
    & $ \argmin_r \max \left| \nu - \tilde{\nu} \right| $
  & 0.21008 & -- \\[7 pt]

    exact at $k = 1$    
    &  {
    $ \cot \left( \frac{\pi}{2} \cdot \frac{1 + e^{-\gamma} - 2 \log 2}{e^{-\gamma} - \log 2 + \frac{1}{3} } \right)  $} 
  & 0.209257 & -- \\[7 pt]

    tangent at $k \approx 0.4184$    
    &  ? 
  & 0.205282 & L \\[4pt]

  \hline

  \end{tabular}
  \label{tab:compareInterps}
\end{table}

\section{Conclusions}
\label{sec:conclusions}

Tight upper and lower bounds to the median of the gamma distribution are introduced.  The simplest new lower bound is never below the 48th percentile, and the simplest new upper bound, of the same form $2^{-1/k}(A + Bk)$, is never above the 55th percentile, over the entire range of $k > 0$.  Using arctan and rational-function interpolators between these bounds, two one-parameter families of closed-form bounds and approximations to the median of a gamma distribution are proposed.

The one-parameter rational-function family has simple closed-form formulae for tightest upper and lower bounds, staying below 50.85 and above 49.69 percentile, respectively; higher-order rational functions can provide tighter bounds or better approximations.  

The one-parameter arctan family of interpolators is a better fit to the ideal interpolator, and includes a version that is most accurate in the low-$k$ tail and provides a closed-from tight upper bound, staying below 50.18 percentile.  With different $b$ parameter, several approximations in the family, including the closed-form version that is exact at $k=1$, stay between 49.97 and 50.03 percentile.  We have not found an analytic formula for the parameter that gives the tight lower bound, which stays above 49.96 percentile, but have shown where to find it, graphically or numerically.

The approach of interpolating between tight bounds opens the way to finding tighter bounds and more accurate approximations, and to finding more such families of bounds and approximations via other interpolator forms.

\section*{Acknowledgements}
The author acknowledges the helpful comments of Google colleagues, especially Pascal Getreuer, Srinivas Vasudevan, and Dan Piponi.

\bibliographystyle{abbrv}

\bibliography{references.bib} 

\begin{thebibliography}{10}

\bibitem{adell2008ramanujan}
J.~Adell and P.~Jodr{\'a}.
\newblock On a {Ramanujan} equation connected with the median of the gamma
  distribution.
\newblock {\em Transactions of the American Mathematical Society},
  360(7):3631--3644, 2008.

\bibitem{barry2000approximation}
D.~Barry, J.-Y. Parlange, and L.~Li.
\newblock Approximation for the exponential integral ({Theis} well function).
\newblock {\em Journal of Hydrology}, 227(1-4):287--291, 2000.

\bibitem{berg2006chen}
C.~Berg and H.~L. Pedersen.
\newblock The {Chen--Rubin} conjecture in a continuous setting.
\newblock {\em Methods and Applications of Analysis}, 13(1):63--88, 2006.

\bibitem{berg2008convexity}
C.~Berg and H.~L. Pedersen.
\newblock Convexity of the median in the gamma distribution.
\newblock {\em Arkiv f{\"o}r Matematik}, 46(1):1--6, 2008.

\bibitem{chen2017median}
C.-P. Chen.
\newblock The median of gamma distribution and a related {Ramanujan} sequence.
\newblock {\em The Ramanujan Journal}, 44(1):75--88, 2017.

\bibitem{chen1986bounds}
J.~Chen and H.~Rubin.
\newblock Bounds for the difference between median and mean of gamma and
  {Poisson} distributions.
\newblock {\em Statistics \& Probability Letters}, 4(6):281--283, 1986.

\bibitem{choi1994medians}
K.~P. Choi.
\newblock On the medians of gamma distributions and an equation of {Ramanujan}.
\newblock {\em Proceedings of the American Mathematical Society},
  121(1):245--251, 1994.

\bibitem{marsaglia1986c249}
J.~C. Marsaglia.
\newblock The incomplete gamma function and {Ramanujan’s} rational
  approximation to {$e^x$}.
\newblock {\em Journal of Statistical Computation and Simulation},
  24(2):163--168, 1986.

\bibitem{steinbrecher2008quantile}
G.~Steinbrecher and W.~T. Shaw.
\newblock Quantile mechanics.
\newblock {\em European Journal of Applied Mathematics}, 19(2):87--112, 2008.

\bibitem{you2017approximation}
X.~You.
\newblock Approximation of the median of the gamma distribution.
\newblock {\em Journal of Number Theory}, 174:487--493, 2017.

\end{thebibliography}

\end{document}